\documentclass[twoside             
]{amsart}
\usepackage[centertags]{amsmath}

\usepackage{amsfonts}
\usepackage{amssymb}
\usepackage{amsthm}
\usepackage{eucal}
\usepackage[cmtip,all]{xy}
\usepackage{graphicx}

\newtheorem{thm}[equation]{Theorem}

\newtheorem{lem}[equation]{Lemma}
\newtheorem{prop}[equation]{Proposition}
\theoremstyle{definition}
\newtheorem{defn}[equation]{Definition}
\theoremstyle{remark}
\newtheorem{rem}[equation]{Remark}

\numberwithin{equation}{section}

\newcommand{\abs}[1]{\left\vert#1\right\vert}
\newcommand{\set}[1]{\left\{#1\right\}}

\newcommand{\To}{\longrightarrow}

\newcommand{\A}{\mathcal{A}}

\newcommand{\C}[1]{\mathbf{#1}}
\newcommand{\D}[1]{\mathcal{#1}}

\newcommand{\relnn}{{\textrm{\rm R}}}
\newcommand{\relss}{{\textrm{\rm S}}}
\newcommand{\reln}[1]{{\textrm{\rm(\relnn #1)}}}
\newcommand{\rels}[1]{{\textrm{\rm(\relss #1)}}}

\newcommand{\EZDIAG}[5]{\xymatrix@C+=1.5cm{*+[r]{#1}
\ar@(u,l)_(0.62){\displaystyle #5}[]
\ar@<.7ex>^-{#3}[r]&\ar@<.7ex>^-{#4}[l]#2}}

\def\ni{{nil}}
\def\abb{{ab}}
\newcommand{\sym}[1]{\mathrm{Sym}(#1)}

\def\r{\rightarrow} 
\def\l{\leftarrow} 

\def\into{\rightarrowtail}
\def\onto{\twoheadrightarrow}

\def\F{\mathcal{F}}

\def\sm{\mathrm{Sum}}

\def\st{\stackrel} 

\renewcommand{\ker}{\operatorname{Ker}}

\def\Z{\mathbb{Z}}

\newcommand{\grupo}[1]{\langle #1\rangle}

\begin{document}

\title{On $K_1$ of a Waldhausen category}%
\author{Fernando Muro \and Andrew Tonks}%
\address{Max-Planck-Institut f\"ur Mathematik, Vivatsgasse 7, 53111 Bonn, Germany}
\email{muro@mpim-bonn.mpg.de}
\address{London Metropolitan University, 166--220 Holloway Road, London N7 8DB, UK}
\email{a.tonks@londonmet.ac.uk}

\thanks{The first author was partially supported
by the project MTM2004-01865 and the MEC postdoctoral fellowship
EX2004-0616, and the second by the MEC-FEDER grant MTM2004-03629.}
\subjclass{19B99, 16E20, 18G50, 18G55}
\keywords{$K$-theory, Waldhausen category, Postnikov invariant, stable quadratic module, crossed complex,
categorical group}

\begin{abstract}
We give a simple representation of all elements in $K_1$ of a Waldhausen category and
prove relations between these representatives which hold in $K_1$.
\end{abstract} \maketitle \tableofcontents

\section*{Introduction}

A general notion of $K$-theory, for a category $\C{W}$ with
cofibrations and weak equivalences,
was defined by Waldhausen \cite{akttsI}
as the homotopy groups of the loop space of a certain simplicial category $wS.\C{W}$,
$$K_n\C{W}\;=\;\pi_n\Omega\abs{wS.\C{W}}\;\cong\;\pi_{n+1}\abs{wS.\C{W}},\;\; n\geq0.$$
Waldhausen $K$-theory generalizes the $K$-theory of an exact category $\C{E}$, defined as
the homotopy groups of $\Omega\abs{Q\C{E}}$, where $Q\C{E}$ is a category defined by Quillen.

Gillet and Grayson defined in \cite{lsqc} a simplicial set $G.\C{E}$ which is a model for $\Omega\abs{Q\C{E}}$.
This allows one to compute $K_1\C{E}$ as a fundamental group, $K_1\C{E}=\pi_1\abs{G.\C{E}}$. Using
the standard techniques for computing $\pi_1$ Gillet and Grayson produced algebraic representatives for arbitrary elements in $K_1\C{E}$.
These representatives were simplified by Sherman \cite{k1ac, k1ec} and further simplified by Nenashev \cite{dsespaek1}.
Nenashev's representatives are pairs of short exact sequences on the same objects,
\begin{equation*}
\xymatrix@R=5pt@C=10pt{
A\ar@<.7ex>@{>->}[r]\ar@<-.7ex>@{>->}[r]&B\ar@<.7ex>@{->>}[r]\ar@<-.7ex>@{->>}[r]&C.}
\end{equation*}
Such a pair gives a loop in $|G.\C{E}|$ which corresponds to a $2$-sphere in $\abs{wS.\C{E}}$, obtained by pasting the $2$-simplices associated
to each short exact sequence along their common boundary.
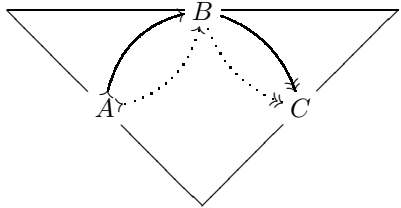
\begin{figure}[h]
\begin{eqnarray*}
\entrymodifiers={+0}
\xymatrix@!=50pt{\ar@{-}[rr]|*+{B}="b"\ar@{-}[rd]|*+{A}="a"&&\\
&\ar@{-}[ru]|*+{C}="c"&
\ar@{>->}@/^10pt/"a";"b"
\ar@{>.>}@/_10pt/"a";"b"
\ar@{->>}@/^8pt/"b";"c"
\ar@{.>>}@/_8pt/"b";"c"}
\end{eqnarray*}
\caption{Nenashev's representative of an element in $K_1\C{E}$.}
\end{figure}

Nenashev showed in \cite{dsesk1} certain relations which are satisfied by these pairs of short exact sequences, and
he later proved in \cite{K1gr} that pairs of short exact sequences together with these relations yield a
presentation of $K_1\C{E}$.

In this paper we produce representatives for all elements in $K_1$ of an arbitrary Waldhausen category $\C{W}$
which are as close as possible to Nenashev's representatives for exact categories. They are given by pairs of
cofiber sequences where the cofibrations have the same source and target. The cofibers, however, may be
non-isomorphic, but they are weakly equivalent via a length $2$ zig-zag.
\begin{equation*}
\xymatrix@R=5pt@C=10pt{&&C_{1}\ar@{<-}[rd]^\sim&\\
A\ar@<.7ex>@{>->}[r]\ar@<-.7ex>@{>->}[r]&B\ar@{->>}[ru]\ar@{->>}[rd]&&C.\\
&&C_{2}\ar@{<-}[ru]_\sim&}
\end{equation*}
This diagram is called a \emph{pair of weak cofiber sequences}. It corresponds to a $2$-sphere in $\abs{wS.\C{W}}$
obtained by pasting the $2$-simplices associated
to the cofiber sequences along the common part of the boundary. The two edges which remain free after this operation
are filled with a disk made of two pieces corresponding to the weak equivalences.
\begin{figure}[h]
\begin{eqnarray*}
\entrymodifiers={+0}
\xymatrix@!=50pt{\ar@{-}[rr]|*+{B}="b"\ar@{-}[rd]|*+{A}="a"&&\\
&\ar@{-}@/^23pt/[ru]|*+{C_1}="c1"\ar@{-}@/_30pt/[ru]|*+{C_2}="c2"\ar@{-}@/_8pt/[ru]|<(.6)*+{C}="c"&
\ar@{>->}@/^10pt/"a";"b"
\ar@{>.>}@/_10pt/"a";"b"
\ar@{->>}@/^2pt/"b";"c1"
\ar@{.>>}@/_8pt/"b";"c2"
\ar@/_3pt/"c";"c1"_-\sim
\ar@/^2pt/"c";"c2"^-\sim}
\end{eqnarray*}
\caption{Our representative of an element in $K_1\C{W}$.}
\end{figure}
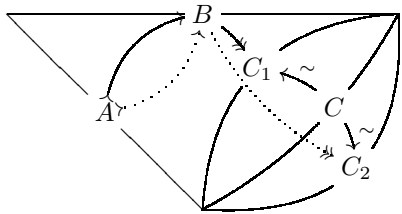

In addition we show in Section \ref{relations} that pairs of weak cofiber sequences satisfy a certain generalization of Nenashev's
relations.

Let $K^\mathrm{wcs}_1\; \C{W}$ be the group generated by the pairs of weak cofiber sequences in $\C{W}$ modulo the relations
in Section \ref{relations}. The results of this paper show in particular that there is a natural surjection
$$K^\mathrm{wcs}_1\; \C{W}\twoheadrightarrow K_1\C{W}$$
which is an isomorphism in case $\C{W}=\C{E}$ is an exact category. We do not know whether it is an
isomorphism for an arbitrary Waldhausen category $\C{W}$, but we prove here further evidence which supports the conjecture:
the homomorphism $$\cdot\,\eta\colon K_0\C{W}\otimes\Z/2\r K_1\C{W}$$ induced by the action of the Hopf map $\eta\in\pi_*^s$
in the stable homotopy groups of spheres  factors as
$$K_0\C{W}\otimes\Z/2\To K^\mathrm{wcs}_1\; \C{W}\twoheadrightarrow K_1\C{W}.$$

We finish this introduction with some comments about the method of proof. We do not rely on any previous
similar result because we do not know of any generalization of the Gillet-Grayson construction for an arbitrary Waldhausen
category $\C{W}$. The largest known class over which the Gillet-Grayson construction works is the class of pseudo-additive categories, which extends the class of exact categories but still does not cover all Waldhausen categories, see \cite{udlakt}. We use instead the algebraic model $\D{D}_*\C{W}$ defined in \cite{1tK}
for the $1$-type of the Waldhausen $K$-theory spectrum $K\C{W}$. The algebraic object $\D{D}_*\C{W}$ is a chain
complex of non-abelian groups concentrated in dimensions $i=0,1$ whose homology is $K_i\C{W}$.
$$\xymatrix{&
(\D{D}_0\C{W})^\abb\otimes
(\D{D}_0\C{W})^\abb\ar[d]_{\grupo{\cdot,\cdot}}&&
\\K_1\C{W}\;\ar@{^{(}->}[r]&
\D{D}_1\C{W}\ar[r]^\partial&\D{D}_0\C{W}\ar@{->>}[r]&K_0\C{W}.}$$
This non-abelian chain complex is equipped with a bilinear map $\grupo{\cdot,\cdot}$ which determines the
commutators. This makes $\D{D}_*\C{W}$ a
stable quadratic module in the sense of \cite{ch4c}. This stable quadratic module is defined in \cite{1tK} in
terms of generators and relations. Generators correspond simply to objects, weak equivalences, and cofiber sequences in
$\C{W}$. The reader will notice that the simplicity of the presentation of $\D{D}_*\C{W}$ makes the proofs of the
results in this paper considerably shorter than the proofs of the analogous results for exact categories.


\subsection*{Acknowledgement}

The authors are very grateful to Grigory Garkusha for asking about the relation between the algebraic model for Waldhausen $K$-theory defined in \cite{1tK} and Nenashev's presentation of $K_1$ of an exact category in \cite{K1gr}.

\section{The algebraic model for $K_0$ and $K_1$}\label{1}

\begin{defn}\label{ob}
A \emph{stable quadratic module} $C_*$ is a diagram of group homomorphisms
$$C_0^\abb\otimes C_0^\abb\st{\grupo{\cdot,\cdot}}\To C_1\st{\partial}\To C_0$$
such that given $c_i,d_i\in C_i$, $i=0,1$,
\begin{enumerate}
\item $\partial\grupo{c_0,d_0}=[d_0,c_0]$,
\item $\grupo{\partial (c_1),\partial (d_1)}=[d_1,c_1]$,
\item $\grupo{c_0,d_0}+\grupo{d_0,c_0}=0$.
\end{enumerate}
Here $[x,y]=-x-y+x+y$ is the commutator of two elements
$x,y\in K$ in any group $K$, and $K^\abb$ is the abelianization of
$K$. It follows from the axioms that the image of $\grupo{\cdot,\cdot}$ and $\ker\partial$ are central in $C_1$, the groups $C_0$ and $C_1$ have nilpotency class $2$, and $\partial(C_1)$ is a normal subgroup of $C_0$.

A \emph{morphism} $f\colon C_*\r D_*$ of stable quadratic modules is given by group homomorphisms $f_i\colon C_i\r D_i$, $i=0,1$, compatible with the structure homomorphisms of $C_*$ and $D_*$, i.e. $f_0\partial=\partial f_1$ and $f_1\grupo{\cdot,\cdot}=\grupo{f_0,f_0}$.
\end{defn}

Stable quadratic modules were introduced in
\cite[Definition IV.C.1]{ch4c}.
Notice, however, that we adopt the opposite
convention for the homomorphism $\grupo{\cdot,\cdot}$.

\begin{rem}\label{crom}
There is a natural right action of $C_0$ on $C_1$ defined by
\begin{eqnarray*}
c_1^{c_0}&=&c_1+\grupo{c_0,\partial(c_1)}.
\end{eqnarray*}
The axioms of a stable quadratic module imply that
commutators in $C_0$ act trivially on $C_1$, and that $C_0$ acts
trivially on the image of $\grupo{\cdot,\cdot}$ and on $\ker\partial$.

The action gives  $\partial:C_1\to C_0$ the structure of a crossed
module. Indeed a stable quadratic module is the same as a commutative
monoid in the category of crossed modules such that the monoid product of two
elements in $C_0$ vanishes when one of them is a commutator, see \cite[Lemma 4.15]{1tK}.

A stable quadratic module can be defined by generators and relations in degree 0 and 1.
In the appendix we give details about the construction of a stable quadratic module defined
by generators and relations. This is useful to understand Definition \ref{LA}
below from a purely group-theoretic perspective.
\end{rem}

We assume the reader has certain familiarity with Waldhausen categories and related concepts. We refer to \cite{wiak} for the basics, see also \cite{akts}.
The following definition was introduced in \cite{1tK}.

\begin{defn}\label{LA}
Let $\C{C}$ be a Waldhausen category with distinguished zero object
$*$, and cofibrations and weak equivalences denoted by
$\into$ and $\st{\sim}\rightarrow$, respectively. A
generic cofiber sequence is denoted by
$$A\into B \onto B/A.$$
We define $\D{D}_*\C{C}$ as the stable quadratic module generated in
dimension zero by the symbols
\begin{itemize}
\item $[A]$ for any object in $\C{C}$,
\end{itemize}
and in dimension one by
\begin{itemize}
\item $[A\st{\sim}\r A']$ for any weak equivalence,
\item $[A\into B\onto B/A]$ for any cofiber
sequence,
\end{itemize}
such that the following relations hold.
\begin{enumerate}
\renewcommand{\labelenumi}{\reln{\arabic{enumi}}}
\item $\partial[A\st{\sim}\r A']=-[A']+[A]$.
\item $\partial[A\into B\onto B/A]=-[B]+[B/A]+[A]$.
\item $[*]=0$.
\item $[A\st{1_A}\r A]=0$.
\item $[A\st{1_A}\r A \onto *]=0$, $[*\into A\st{1_A}\r
A]=0$.
\item For any pair of composable weak equivalences $ A\st{\sim}\r B\st{\sim}\r C$,
$$[A\st{\sim}\r C]=[B\st{\sim}\r C]+[A\st{\sim}\r B].$$
\item For any commutative diagram in $\C{C}$ as follows
$$\xymatrix{A\;\ar@{>->}[r]\ar[d]^\sim&B\ar@{->>}[r]\ar[d]^\sim&B/A\ar[d]^\sim\\A'\;\ar@{>->}[r]&B'\ar@{->>}[r]&B'/A'}$$
we have
\begin{eqnarray*}
[A\st{\sim}\r A']+[B/A\st{\sim}\r B'/A']^{[A]}&=&
-[A'\into
B'\onto B'/A']\\&&+[B\st{\sim}\r B']\\&&+[A\into B\onto B/A].
\end{eqnarray*}
\item For any commutative diagram consisting of four cofiber sequences in $\C{C}$ as follows
$$\xymatrix{&&C/B\\&B/A\;\ar@{>->}[r]&C/A\ar@{->>}[u]\\A\;\ar@{>->}[r]&B\;\ar@{>->}[r]\ar@{->>}[u]&C\ar@{->>}[u]}$$
we have
\begin{eqnarray*}
[B\into C\onto C/B]&&\\+[A\into
B\onto B/A]&=&[A\into C\onto
C/A]\\
&&+[B/A\into C/A\onto C/B]^{[A]}.
\end{eqnarray*}
\item For any pair of objects $A, B$ in $\C{C}$
$$\grupo{[A],[B]}=
-[B\st{i_2}\into A\vee B\st{p_1}\onto A]
+[A\st{i_1}\into A\vee B\st{p_2}\onto B]
.$$
Here $$\xymatrix{A\ar@<.5ex>[r]^-{i_1}&A\vee
B\ar@<.5ex>[l]^-{p_1}\ar@<-.5ex>[r]_-{p_2}&B\ar@<-.5ex>[l]_-{i_2}}$$
are the inclusions and projections of a
coproduct
in $\C{C}$.
\end{enumerate}

\begin{rem}
These relations are quite natural, and some illustration of their meaning
is given in \cite[Figure 2]{1tK}.
They are not however minimal: relation $\reln3$ follows from
$\reln2$ and $\reln5$, and $\reln4$ follows from $\reln6$.
Also $\reln5$ is equivalent to
\begin{itemize}
\item[\reln{$5'$}] $[{*}\into {*}\onto{*}]=0$.
\end{itemize}
This equivalence follows from $\reln8$ on considering the diagrams
$$\xymatrix{&&{*}\\&{*}\;\ar@{>->}[r]&{*}\ar@{->>}[u]\\
A\;\ar@{>->}[r]^{1_A}&A\;\ar@{>->}[r]^{1_A}\ar@{->>}[u]&A\ar@{->>}[u]}\qquad\qquad
\xymatrix{&&{A}\\&{*}\;\ar@{>->}[r]&{A}\ar@{->>}[u]_{1_A}\\
{*}\;\ar@{>->}[r]&{*}\;\ar@{>->}[r]\ar@{->>}[u]&A\ar@{->>}[u]_{1_A}}
$$
\end{rem}

We now introduce the basic object of study of this paper.
\begin{defn}
A \emph{weak cofiber sequence} in a Waldhausen category $\C{C}$ is a diagram
$$A\into B\onto C_{1}\st{\sim}\l C$$
given by a cofiber sequence followed by a weak equivalence in the opposite direction.
\end{defn}
We associate the element
\begin{eqnarray}\label{weakcofiber}
[A\into B\onto C_{1}\st{\sim}\l C]&=&[A\into B\onto C_{1}]+[C\st{\sim}\r C_{1}]^{[A]} \in\D{D}_*\C{C}
\end{eqnarray}
to any weak cofiber sequence.
By \reln1 and \reln2 we have
\begin{eqnarray*}
\partial [A\into B\onto C_{1}\st{\sim}\l C]&=&-[B]+[C]+[A].
\end{eqnarray*}

The following fundamental identity for weak cofiber sequences is a slightly less trivial consequence of the relations \reln1--\reln9.
\begin{prop}\label{weaknena}
Consider a commutative diagram
\[\xymatrix{
 A'  \ar@{>->}[r]^-{j^A}
     \ar@{>->}[d]_-{j'}
&A   \ar@{->>}[r]^-{r^A}
     \ar@{>->}[d]_-{j}
&A'''\ar@{>->}[d]_-{j'''}
&A'' \ar@{>->}[d]^-{j''}
     \ar      [l]_-\sim^-{w^A}
\\
 B'  \ar@{>->}[r]^-{j^B}
     \ar@{->>}[d]_-{r'}
&B   \ar@{->>}[r]^-{r^B}
     \ar@{->>}[d]_-{r}
&B'''\ar@{->>}[d]_-{r'''}
&B'' \ar@{->>}[d]^-{r''}
     \ar      [l]_-\sim^-{w^B}
\\
 D'  \ar@{>->}[r]^-{j^D}
&D   \ar@{->>}[r]^-{r^D}
&D'''
&D'' \ar      [l]_-\sim^-{w^D}
\\
 C'  \ar@{>->}[r]^-{j^C}
     \ar      [u]_-\sim^-{w'}
&C   \ar@{->>}[r]^-{r^C}
     \ar      [u]_-\sim^-{w}
&C'''\ar      [u]_-\sim^-{w'''}
&C'' \ar      [u]_-\sim^-{w''}
     \ar      [l]_-\sim^-{w^C}
}
\]
Denote the horizontal weak cofiber sequences by $l^A$, $l^B$, $l^D$ and $l^C$,
and the vertical ones by $l'$, $l$, $l'''$ and $l''$.
Then the following equation holds in $\D{D}_*\C{C}$.
\begin{eqnarray*}
-[l^A]-[l^C]^{[A]}-[l]+[l^B]+[l'']^{[B']}+[l']&=&\grupo{[C'],[A'']}.
\end{eqnarray*}
\end{prop}

\begin{proof}
Applying \reln7 to the two obvious equivalences of cofiber sequences gives
\begin{align*}
\tag{a}[D'\into D\onto D''']&=[w]+[l^C]-[w^C]^{[C']}-[w''']^{[C']}-[w'],\\
\tag{b}[A'''\into B'''\onto D''']&=[w^B]+[l'']-[w'']^{[A'']}-[w^D]^{[A'']}-[w^A],
\end{align*}
using the notation of \eqref{weakcofiber}.

Now consider the following commutative diagrams
\begin{equation*}
\begin{array}{ccc}
\xymatrix{A'\ar@{>->}[r]^{j^A}\ar@{>->}[d]_{j'}\ar@{}[rd]|{\text{push}}
&A\ar@{>->}[d]\ar@/^15pt/@{>->}[rdd]^j&\\
B'\ar@{>->}[r]\ar@/_15pt/@{>->}[rrd]_{j^B}&X\ar@{>->}[rd]&\\
&&B\,,}
&\qquad&
\xymatrix{A'\ar@{>->}[r]^{j^A}\ar@{>->}[d]_{j'}\ar@{}[rd]|{\text{push}}
&A\ar@{>->}[d]\ar@/^15pt/[rdd]^{i_1r^A}&\\
B'\ar@{>->}[r]\ar@/_15pt/[rrd]_{i_2r'}&X\ar@{->>}[rd]&\\
&&A'''\vee D'\,.}
\end{array}
\end{equation*}
We now have four commutative diagrams of cofiber sequences as in \reln8.
\begin{equation*}
\begin{array}{ccccc}
\textrm{c)}
&\xymatrix{&&D'\\
&A'''\;\ar@{>->}[r]&A'''\vee D'\ar@{->>}[u]\\
A'\;\ar@{>->}[r]&A\;\ar@{>->}[r]\ar@{->>}[u]&X\ar@{->>}[u]}
&\;\;\;\;\;\;\;\;\;&
\textrm{d)}
&\xymatrix{&&A'''\\
&D'\;\ar@{>->}[r]&A'''\vee D'\ar@{->>}[u]\\
  A'\;\ar@{>->}[r]&B'\;\ar@{>->}[r]\ar@{->>}[u]&X\ar@{->>}[u]}\\[1cm] \\
\textrm{e)}
&\xymatrix{&&D'''\\
&D'\;\ar@{>->}[r]&D\ar@{->>}[u]\\
A\;\ar@{>->}[r]&X\;\ar@{>->}[r]\ar@{->>}[u]&B\ar@{->>}[u]}
&\;\;\;\;\;\;\;\;\;&
\textrm{f)}
&\xymatrix{&&D'''\\
&A'''\;\ar@{>->}[r]&B'''\ar@{->>}[u]\\B'\;\ar@{>->}[r]&X\;\ar@{>->}[r]\ar@{->>}[u]&B\ar@{->>}[u]}
\end{array}
\end{equation*}
Therefore we obtain the corresponding relations
\begin{align*}
\tag{c}[A\into X\onto D']&+[l^A]-[w^A]^{[A']}
\\&=[A'\into X\onto A'''\vee D']+[A'''\into A'''\vee D'\onto D']^{[A']},
\\\tag{d}[B'\into X\onto A''']&+[l']-[w']^{[A']}
\\&=[A'\into X\onto A'''\vee D']+[D'\into A'''\vee D'\onto A''']^{[A']},
\\\tag{e}[X\into B\onto D''']&+[A\into X\onto D']
\\&=[l]-[w]^{[A]}+[D'\into D\onto D''']^{[A]}
\\&=[l]+([l^C]-[w^C]^{[C']}-[w''']^{[C']}-[w'])^{[A]},
\\\tag{f}[X\into B\onto D''']&+[B'\into X\onto A''']
\\&=[l^B]-[w^B]^{[B']}+[A'''\into B'''\onto D''']^{[B']}
\\&=[l^B]+([l'']-[w'']^{[A'']}-[w^D]^{[A'']}-[w^A])^{[B']}.
\end{align*}
Here we have used equations (a) and (b).
Now $-({\rm c})+({\rm d})$ and $-({\rm e})+({\rm f})$ yield
\begin{align*}
\tag{g}[w^A]^{[A']}-[l^A]
-[A\into X\onto D']\;+\;&[B'\into X\onto A''']
+[l']-[w']^{[A']}=\\
{}-[A'''\into A'''\vee D'\onto D']^{[A']}\;+\;&[D'\into A'''\vee D'\onto A''']^{[A']}=\grupo{[D'],[A''']},
\\[1.5ex]
\tag{h}{}-[A\into X\onto D']\;+\;&[B'\into X\onto A''']=\\
[w']^{[A]}+a^{[C']+[A]}
-[l^C]^{[A]}-[l]\;+\;&[l^B]+[l'']^{[B']}
-a^{[A'']+[B']}
-[w^A]^{[B']},
\end{align*}
where we use \reln9 and Remark \ref{crom}, and we write $a=[w''']+[w^C]=[w^D]+[w'']$, by \reln6.
In fact it is easy to check by using Remark \ref{crom} and the laws of stable quadratic modules that the terms $a^{[C']+[A]}$ and $a^{[A'']+[B']}$ cancel in this expression, since
\[
\partial\left({}-[l^C]^{[A]}-[l]\;+\;[l^B]+[l'']^{[B']}\right)\;\;=\;\;-([C']+[A])+([A'']+[B']).
\]
Substituting (h) into the left hand side of (g) we obtain
\begin{align*}
&[w^A]^{[A']}-[l^A]
+[w']^{[A]}
-[l^C]^{[A]}-[l]+[l^B]+[l'']^{[B']}
-[w^A]^{[B']}
+[l']-[w']^{[A']}
\\&\qquad={}\grupo{[D'],[A''']},
\end{align*}
which, using again Definition \ref{ob} and Remark \ref{crom}, may be rewritten as
\begin{align*}
&
-[l^A]
-[l^C]^{[A]}-[l]+[l^B]+[l'']^{[B']}
+[l']
\\&\qquad={}
-[w']^{[A'']+[A']}-[w^A]^{[A']}
-\grupo{[A'''],[D']}
+[w']^{[A']}+[w^A]^{[C']+[A']}
.
\end{align*}
The result then follows from the identity
\[
\grupo{\partial[w^A],\partial[w']}-
\grupo{[A''],\partial[w']}+
\grupo{[C'],\partial[w^A]}
=\grupo{[C'],[A'']}+
\grupo{[A'''],[D']}.
\]
\end{proof}

\section{Pairs of weak cofiber sequences and $K_1$}

Let $\C{C}$ be a Waldhausen category. A \emph{pair of weak cofiber sequences} is a diagram given by two weak cofiber sequences with the first, second, and fourth objects in common
\begin{equation*}
\xymatrix@R=5pt@C=10pt{&&C_{1}\ar@{<-}[rd]^\sim&\\
A\ar@<.7ex>@{>->}[r]\ar@<-.7ex>@{>->}[r]&B\ar@{->>}[ru]\ar@{->>}[rd]&&C.\\
&&C_{2}\ar@{<-}[ru]_\sim&}
\end{equation*}
We associate to any such pair the following element
\begin{eqnarray}
\label{pwcs}&&\\ \nonumber \set{\begin{array}{c}\xymatrix@R=5pt@C=10pt{&&C_{1}\ar@{<-}[rd]^\sim_{w_1}&\\
A\ar@<.7ex>@{>->}[r]^{j_1}\ar@<-.7ex>@{>->}[r]_{j_2}&B\ar@{->>}[ru]^{r_1}\ar@{->>}[rd]_{r_2}&&C\\
&&C_{2}\ar@{<-}[ru]_\sim^{w_2}&}\end{array}}&=&-[A\st{j_1}\into B\st{r_1}\onto C_{1}\mathop{\l}\limits^\sim_{w_1} C]+[A\st{j_2}\into B\st{r_2}\onto C_{2}\mathop{\l}\limits^\sim_{w_2} C]\\
\nonumber &=&[A\st{j_2}\into B\st{r_2}\onto C_{2}\mathop{\l}\limits^\sim_{w_2} C]-[A\st{j_1}\into B\st{r_1}\onto C_{1}\mathop{\l}\limits^\sim_{w_1} C],
\end{eqnarray}
which lies in $K_1\C{C}=\ker\partial$. For the second equality in (\ref{pwcs}) we use that $\ker\partial$ is central, see Remark \ref{crom}, so we can permute the terms cyclically. 

The following theorem is one of the main results of this paper.

\begin{thm}\label{main}
Any element in $K_1\C{C}$ is represented by a pair of weak cofiber sequences.
\end{thm}

This theorem will be proved later in this paper. We first give a set of useful relations between pairs of weak cofiber sequences and develop a sum-normalized version of the model $\D{D}_*\C{C}$.

\section{Relations between pairs of weak cofiber sequences}\label{relations}

Suppose that we have six pairs of weak cofiber sequences in a Waldhausen category $\C{C}$
\begin{equation*}
\begin{array}{ccc}
\xymatrix@R=5pt@C=10pt{&&A''_{1}\ar@{<-}[rd]_\sim^{w^A_1}&\\
A'\ar@<.7ex>@{>->}[r]^{j^A_1}\ar@<-.7ex>@{>->}[r]_{j^A_2}&A\ar@{->>}[ru]^{r^A_1}\ar@{->>}[rd]_{r^A_2}&&A'',\\
&&A''_{2}\ar@{<-}[ru]^\sim_{w^A_2}&}
&
\xymatrix@R=5pt@C=10pt{&&B''_{1}\ar@{<-}[rd]_\sim^{w^B_1}&\\
B'\ar@<.7ex>@{>->}[r]^{j^B_1}\ar@<-.7ex>@{>->}[r]_{j^B_2}&B\ar@{->>}[ru]^{r^B_1}\ar@{->>}[rd]_{r^B_2}&&B'',\\
&&B''_{2}\ar@{<-}[ru]^\sim_{w^B_2}&}
&
\xymatrix@R=5pt@C=10pt{&&\check{C}''_{1}\ar@{<-}[rd]_\sim^{w^C_1}&\\
C'\ar@<.7ex>@{>->}[r]^{j^C_1}\ar@<-.7ex>@{>->}[r]_{j^C_2}&C\ar@{->>}[ru]^{r^C_1}\ar@{->>}[rd]_{r^C_2}&&C'',\\
&&\check{C}''_{2}\ar@{<-}[ru]^\sim_{w^C_2}&}
\\
\xymatrix@R=5pt@C=10pt{&&C'_{1}\ar@{<-}[rd]_\sim^{w'_1}&\\
A'\ar@<.7ex>@{>->}[r]^{j'_1}\ar@<-.7ex>@{>->}[r]_{j'_2}&B'\ar@{->>}[ru]^{r'_1}\ar@{->>}[rd]_{r'_2}&&C',\\
&&C'_{2}\ar@{<-}[ru]^\sim_{w'_2}&}
&
\xymatrix@R=5pt@C=10pt{&&C_{1}\ar@{<-}[rd]_\sim^{w_1}&\\
A\ar@<.7ex>@{>->}[r]^{j_1}\ar@<-.7ex>@{>->}[r]_{j_2}&B\ar@{->>}[ru]^{r_1}\ar@{->>}[rd]_{r_2}&&C,\\
&&C_{2}\ar@{<-}[ru]^\sim_{w_2}&}
&
\xymatrix@R=5pt@C=10pt{&&\hat{C}''_{1}\ar@{<-}[rd]_\sim^{w''_1}&\\
A''\ar@<.7ex>@{>->}[r]^{j''_1}\ar@<-.7ex>@{>->}[r]_{j''_2}&B''\ar@{->>}[ru]^{r''_1}\ar@{->>}[rd]_{r''_2}&&C'',\\
&&\hat{C}''_{2}\ar@{<-}[ru]^\sim_{w''_2}&}
\end{array}
\end{equation*}
that we denote for the sake of simplicity as
$\lambda^A$, $\lambda^B$, $\lambda^C$, $\lambda'$, $\lambda$, and $\lambda''$, respectively.
Moreover, assume that for $i=1,2$
there exists a commutative diagram
\begin{equation*}
\xymatrix{
A'\ar@{>->}[r]^-{j^A_i}\ar@{>->}[d]_-{j'_i}&
A\ar@{>->}[d]_-{j_i}\ar@{->>}[r]^-{r^A_{i}}&A''_i\ar@{>->}[d]_-{\check\jmath_{i}}
&A''\ar@{>->}[d]^-{j''_i}\ar[l]_-\sim^-{w^A_{i}}
\\
B'\ar@{>->}[r]^-{j^B_i}\ar@{->>}[d]_-{r'_{i}}&B\ar@{->>}[r]^-{r^B_{i}}\ar@{->>}[d]_-{r_{i}}&
B''_i\ar@{->>}[d]_-{\check r_{i}}
&B''\ar@{->>}[d]^-{r''_i}\ar[l]_-\sim^-{w^B_{i}}
\\
C'_i\ar@{>->}[r]^-{\hat\jmath_{i}}&C_i\ar@{->>}[r]^-{\hat r_{i}}&C''_i
&\hat{C}''_i\ar[l]_-\sim^-{\hat{w}_i}
\\
C'\ar@{>->}[r]_-{j^C_i}\ar[u]_-\sim^{w'_{i}}&C\ar@{->>}[r]_-{r^C_i}\ar[u]_-\sim^{w_{i}}&
\check{C}''_i\ar[u]_-\sim^{\check{w}_i}
&C''\ar[l]_-\sim^-{w^C_{i}}\ar[u]_-\sim^{w''_{i}}}
\end{equation*}
Notice that six of the eight weak cofiber sequences in this diagram are already in the six diagrams above. The other two weak cofiber sequences
are just assumed to exist with the property of making the diagram commutative.

\begin{thm}\label{larel}
In the situation above the following equation holds
\begin{align*}\tag{\relss1}
\set{\lambda^A}-\set{\lambda^B}+\set{\lambda^C}&=\set{\lambda'}-\set{\lambda}+\set{\lambda''}.
\end{align*}
\end{thm}
\begin{proof}
For $i=1,2$ let $\lambda^A_i$, $\lambda^B_i$, $\lambda^C_i$, $\lambda'_i$, $\lambda_i$, and $\lambda''_i$ be the upper and the lower weak cofiber sequences of the pairs above, respectively. By Proposition \ref{weaknena}
$$\begin{array}{l}
-[\lambda_1^A]-[\lambda_1^C]^{[A]}-[\lambda_1]+[\lambda_1^B]+[\lambda_1'']^{[B']}+[\lambda_1']\;\;=\;\;\grupo{[C'],[A'']}\\
\qquad\qquad\qquad\qquad\qquad\qquad=\;\;-[\lambda_2^A]-[\lambda_2^C]^{[A]}-[\lambda_2]+[\lambda_2^B]+[\lambda_2'']^{[B']}+[\lambda_2'].
\end{array}$$
Since $\partial(\lambda^C_1)=\partial(\lambda^C_2)$ and $\partial(\lambda''_1)=\partial(\lambda''_2)$ the actions cancel,
$$-[\lambda_1^A]-[\lambda_1^C]-[\lambda_1]+[\lambda_1^B]+[\lambda_1'']+[\lambda_1']
\;\;=\;\;-[\lambda_2^A]-[\lambda_2^C]-[\lambda_2]+[\lambda_2^B]+[\lambda_2'']+[\lambda_2'].$$
Now one can rearrange the terms in this equation, by using that pairs of weak cofiber sequences are central in $\D{D}_1\C{C}$, obtaining the equation in the statement.
\end{proof}

Theorem \ref{larel} encodes the most relevant relation satisfied by pairs of weak cofiber sequences in $K_1\C C$. They satisfy a further relation which follows from the very definition in (\ref{pwcs}), but which we would like to record as a proposition by analogy with \cite{dsesk1}.

\begin{prop}\label{las0}
A pair of weak cofiber sequences given by two times the same weak cofiber sequence is trivial.
\begin{align*}\tag{\relss2}
\set{\begin{array}{c}\xymatrix@R=5pt@C=10pt{&&C_{1}\ar@{<-}[rd]^\sim_{w}&\\
A\ar@<.7ex>@{>->}[r]^{j}\ar@<-.7ex>@{>->}[r]_{j}&B\ar@{->>}[ru]^{r}\ar@{->>}[rd]_{r}&&C\\
&&C_{1}\ar@{<-}[ru]_\sim^{w}&}\end{array}}&=0.
\end{align*}
\end{prop}

We establish a further useful relation in the following proposition.

\begin{prop}\label{spwcs}
Relations \rels1 and \rels2 imply that the sum of two pairs of weak cofiber sequences coincides with their coproduct.
\begin{eqnarray*}
\left\{\begin{array}{c}\xymatrix@R=5pt@C=10pt{&&C_{1}\vee\bar{C}_{1}\ar@{<-}[rd]^\sim&\\
A\vee\bar{A}\ar@<.7ex>@{>->}[r]\ar@<-.7ex>@{>->}[r]&B\vee\bar{B}\ar@{->>}[ru]\ar@{->>}[rd]&&C\vee\bar{C}\\
&&C_{2}\vee\bar{C}_{2}\ar@{<-}[ru]_\sim&}\end{array}\right\}&=&\left\{\begin{array}{c}\xymatrix@R=5pt@C=10pt{&&C_{1}\ar@{<-}[rd]^\sim&\\
A\ar@<.7ex>@{>->}[r]\ar@<-.7ex>@{>->}[r]&B\ar@{->>}[ru]\ar@{->>}[rd]&&C\\
&&C_{2}\ar@{<-}[ru]_\sim&}\end{array}\right\}\\
&&+\left\{\begin{array}{c}\xymatrix@R=5pt@C=10pt{&&\bar{C}_{1}\ar@{<-}[rd]^\sim&\\
\bar{A}\ar@<.7ex>@{>->}[r]\ar@<-.7ex>@{>->}[r]&\bar{B}\ar@{->>}[ru]\ar@{->>}[rd]&&\bar{C}\\
&&\bar{C}_{2}\ar@{<-}[ru]_\sim&}\end{array}\right\}.
\end{eqnarray*}
\end{prop}

\begin{proof}
Apply relations \rels1 and \rels2 to the following pairs of weak cofiber sequences
\begin{equation*}
\begin{array}{ccc}
\!\!\!\!\!\!\!\!\!\!\!\!
\xymatrix@R=5pt@C=10pt{&&C_{1}\ar@{<-}[rd]^\sim&\\
A\ar@<.7ex>@{>->}[r]\ar@<-.7ex>@{>->}[r]&B\ar@{->>}[ru]\ar@{->>}[rd]&&C,\\
&&C_{2}\ar@{<-}[ru]_\sim&}
&\!\!\!\!\!\!\!\!\!\!\!\!
\xymatrix@R=5pt@C=10pt{&&C_{1}\vee\bar{C}_{1}\ar@{<-}[rd]^\sim&\\
A\vee\bar{A}\ar@<.7ex>@{>->}[r]\ar@<-.7ex>@{>->}[r]&B\vee\bar{B}\ar@{->>}[ru]\ar@{->>}[rd]&&C\vee\bar{C},\\
&&C_{2}\vee\bar{C}_{2}\ar@{<-}[ru]_\sim&}
&\!\!\!\!\!
\xymatrix@R=5pt@C=10pt{&&\bar{C}_{1}\ar@{<-}[rd]^\sim&\\
\bar{A}\ar@<.7ex>@{>->}[r]\ar@<-.7ex>@{>->}[r]&\bar{B}\ar@{->>}[ru]\ar@{->>}[rd]&&\bar{C},\\
&&\bar{C}_{2}\ar@{<-}[ru]_\sim&}
\\
\!\!\!\!
\xymatrix@R=5pt@C=10pt{&&\bar{A}\ar@{<-}[rd]_-\sim^-{1}&\\
A\ar@<.7ex>@{>->}[r]^-{i_1}\ar@<-.7ex>@{>->}[r]_-{i_1}&A\vee\bar{A}\ar@{->>}[ru]^-{p_2}\ar@{->>}[rd]_-{p_2}&&\bar{A},\\
&&\bar{A}\ar@{<-}[ru]^-\sim_-{1}&}
&\!\!\!\!\!\!\!\!\!\!\!\!\!\!\!\!
\xymatrix@R=5pt@C=10pt{&&\bar{B}\ar@{<-}[rd]_-\sim^-{1}&\\
B\ar@<.7ex>@{>->}[r]^-{i_1}\ar@<-.7ex>@{>->}[r]_-{i_1}&B\vee\bar{B}\ar@{->>}[ru]^-{p_2}\ar@{->>}[rd]_-{p_2}&&\bar{B},\\
&&\bar{B}\ar@{<-}[ru]^-\sim_-{1}&}
&\!\!\!\!\!\!\!\!\!\!\!\!\!\!\!\!
\xymatrix@R=5pt@C=10pt{&&\bar{C}\ar@{<-}[rd]_-\sim^-{1}&\\
C\ar@<.7ex>@{>->}[r]^-{i_1}\ar@<-.7ex>@{>->}[r]_-{i_1}&C\vee\bar{C}\ar@{->>}[ru]^-{p_2}\ar@{->>}[rd]_-{p_2}&&\bar{C}.\\
&&\bar{C}\ar@{<-}[ru]^-\sim_-{1}&}
\end{array}
\end{equation*}
\end{proof}

\section{Waldhausen categories with functorial coproducts}

Let $\C{C}$ be a Waldhausen category endowed with a symmetric monoidal structure $+$
which is strictly associative
\begin{eqnarray*}
(A+B)+C&=&A+(B+C),
\end{eqnarray*}
strictly unital
\begin{eqnarray*}
{*}+A\;\;=\;\;A&=&A+{*},
\end{eqnarray*}
but not necessarily strictly commutative, such that
$$A\;=\;A+{*}\To A+B\longleftarrow {*}+B\;=\;B$$
is always a coproduct diagram. Such a category will be called a
\emph{Waldhausen category with a functorial coproduct}. Then we define
the sum-normalized stable quadratic module $\D{D}_*^+\C{C}$ as the
quotient of $\D{D}_*\C{C}$ by the extra relation
\begin{enumerate}
\item[\reln{10}] $[B\st{i_2}\into A+ B\st{p_1}\onto A]=0$.
\end{enumerate}
\end{defn}

\begin{rem}\label{A+B}
In $\D{D}_*^+\C{C}$, the relations $\reln2$ and $\reln{10}$ imply that
$$
[A+B]=[A]+[B]\,.
$$
Furthermore, relation $\reln9$ becomes equivalent to
\begin{itemize}
\item[\reln{$9'$}]  $\grupo{[A],[B]}=[\tau_{B,A}\colon B+A\st{\cong}\To A+B]$.
\end{itemize}
Here $\tau_{B,A}$ is the symmetry isomorphism of the symmetric monoidal structure.
This equivalence follows from the commutative diagram
$$\xymatrix{A\;\ar@{>->}[r]^-{i_2}\ar@{=}[d]&B+A\ar@{->>}[r]^-{p_1}\ar[d]_{\tau_{B,A}}^\cong&B\ar@{=}[d]\\
A\;\ar@{>->}[r]_-{i_1}&A+B\ar@{->>}[r]_-{p_2}&B}$$
together with $\reln4$, $\reln7$ and $\reln{10}$.
\end{rem}

\begin{thm}\label{sumnor}
Let $\C{C}$ be a Waldhausen category with a functorial coproduct such
that there exists a set $\mathbb{S}$ which freely generates the monoid
of objects of $\C C$ under the operation $+$. Then the projection
$$p\colon \D{D}_*\C{C}\onto \D{D}_*^+\C{C}$$
is a weak equivalence. Indeed, it is part of a strong deformation retraction.
\end{thm}

Waldhausen categories satisfying the hypothesis of Theorem \ref{sumnor} are general enough as the following
proposition shows.

\begin{prop}\label{hay}
For any Waldhausen category $\C{C}$ there is another Waldhausen category $\sm(\C{C})$ with a functorial coproduct
whose monoid of objects is freely generated by the objects of $\C{C}$ except from $*$. Moreover, there are
natural mutually inverse exact equivalences of categories
$$\xymatrix{\sm(\C{C})\ar@<.5ex>[r]^-\varphi&\C{C}\ar@<.5ex>[l]^-\psi}.$$
\end{prop}

\begin{proof}
The statement already says which are the objects of $\sm(\C{C})$. The functor $\varphi$ sends an object $Y$ in $\sm(\C{C})$, which can be uniquely written as a formal sum of non-zero objects in $\C{C}$, $Y=X_1+\cdots+X_n$, to an arbitrarily chosen coproduct of these objects in $\C{C}$
\begin{eqnarray*}
\varphi(Y)&=& X_1\vee\cdots\vee X_n.
\end{eqnarray*}
For $n=0,1$ we make special choices, namely for $n=0$, $\varphi(Y)=*$, and for $n=1$ we set $\varphi(Y)=X_1$.
Morphisms in $\sm(\C{C})$ are defined in the unique possible way so that $\varphi$ is fully faithful. Then the formal sum defines a functorial coproduct on $\sm(\C{C})$. The functor $\psi$ sends $*$ to the zero object of the free monoid and any other object in $\C{C}$ to the corresponding object with a single summand in $\sm(\C{C})$, so that $\varphi\psi=1$. The inverse natural isomorphism $1\cong\psi\varphi$,
\begin{eqnarray*}
X_1+\cdots+ X_n&\cong&(X_1\vee\cdots\vee X_n),
\end{eqnarray*}
is the unique isomorphism in $\sm(\C{C})$ which $\varphi$ maps to the identity. Cofibrations and weak equivalences in $\sm(\C{C})$ are the morphisms which $\varphi$ maps to cofibrations and weak equivalences, respectively. This Waldhausen category structure in $\sm(\C{C})$ makes the functors $\varphi$ and $\psi$ exact.
\end{proof}

Let us recall the notion of homotopy in the category of stable quadratic modules.

\begin{defn}\label{homo}
Given morphisms $f,g:C_*\r C_*'$ of stable quadratic modules,
a \emph{homotopy} $f\st\alpha\simeq g$ from $f$ to $g$ is a function
$\alpha\colon C_0\r C_1'$ satisfying
\begin{enumerate}
\item $\alpha(c_0+d_0)=\alpha(c_0)^{f_0(d_0)}+\alpha(d_0)$,
\item $g_0(c_0)=f_0(c_0)+\partial\alpha(c_0)$,
\item $g_1(c_1)=f_1(c_1)+\alpha\partial(c_1)$.
\end{enumerate}
\end{defn}

The following lemma then follows from the laws of stable quadratic modules.

\begin{lem}\label{-c}
A homotopy $\alpha$ as in Definition \ref{homo} satisfies
\begin{align}
\tag{a}\alpha([c_0,d_0])&=-\grupo{f_0(d_0),f_0(c_0)}+\grupo{g_0(d_0),g_0(c_0)},\\
\tag{b}\alpha(c_0)+g_1(c_1)&=f_1(c_1)+\alpha(c_0+\partial(c_1)).
\end{align}
\end{lem}

Now we are ready to prove Theorem \ref{sumnor}.

\begin{proof}[Proof of Theorem \ref{sumnor}]
In order to define a strong deformation retraction,
$$\EZDIAG{\D{D}_*\C{C}}{\D{D}_*^+\C{C}}pj\alpha\;,
\qquad
1\st\alpha\simeq jp,\;\;\;pj=1,\;\;
$$
the crucial step will be to define the homotopy $\alpha\colon\D{D}_0\C{C}\r\D{D}_1\C{C}$ on the generators $[A]$ of
$\D{D}_0\C{C}$.
Then one can use the equations
\begin{enumerate}
\item $\alpha(c_0+d_0)=\alpha(c_0)^{d_0}+\alpha(d_0)$,
\item $jp(c_0)=c_0+\partial\alpha(c_0)$,
\item $jp(c_1)=c_1+\alpha\partial(c_1)$,
\end{enumerate}
for $c_i,d_i\in\D{D}_i\C{C}$ to define the
composite $jp$ and the homotopy $\alpha$ on all of $\D{D}_*\C{C}$.
It is a straightforward calculation to check that the
map $jp$ so defined is a homomorphism and that $\alpha$ is well
defined. In particular, Lemma \ref{-c} (a) implies that $\alpha$ vanishes on commutators of length 3.
In order to define $j$ we must show that $jp$ factors through the
projection $p$, that
is,
$$
jp([B\st{i_2}\into A+ B\st{p_1}\onto A])\;\;=\;\;0.
$$
If we also show that $$p\alpha\;\;=\;\;0$$
then equations (2) and (3) above say
$pjp(c_i)=p(c_i)+0$, and since $p$ is surjective it follows that the composite $pj$ is
the identity.

The set of objects of $\C C$ is the free monoid generated by objects
$S\in \mathbb S$. Therefore we can define inductively the homotopy $\alpha$ by
$$\alpha([S+B])=
[B\st{i_2}\into S+B\st{p_1}\onto S]
+\alpha([B])$$ for $B\in\C C$ and
$S\in\mathbb S$, with $\alpha([S])=0$.
We claim that a similar relation then holds for all objects $A$, $B$ of $\C C$,
\begin{align}
\label{alphaAB}\alpha([A+B])
&=
[B\st{i_2}\into A+B\st{p_1}\onto A]+\alpha([A]+[B])
\\ \nonumber
&=
[B\st{i_2}\into A+B\st{p_1}\onto A]+
\alpha([A])^{[B]}+\alpha([B])
.
\end{align}
If $A=*$ or $A\in\mathbb S$ then \eqref{alphaAB} holds by definition, so
we assume inductively it holds for given $A$, $B$ and show
it holds for $S+A$, $B$ also:
\begin{align*}
&\alpha([S+A+B])=
[A+B\into S+A+B\onto S]
+\alpha([A+B])\\
&\qquad=
[A+B\into S+A+B\onto S]
+[B\into A+B\onto A]\;\;+
\;\alpha([A])^{[B]}+\alpha([B])\\
&\qquad=
[B\into S+A+B\onto S+A]+
[A\into S+A\onto S]^{[B]}
+\alpha([A])^{[B]}+\alpha([B])\\
&\qquad=
[B\into S+A+B\onto S+A]+
\alpha([S+A])^{[B]}+\alpha([B])\,.
\end{align*}
Here we have used $\reln8$ for the composable
cofibrations
$$
B\;\into\; A+B\;\into\; S+A+B\,.
$$
It is clear from the definition of $\alpha$ that the relation $p\alpha=0$ holds.
It remains to see that $jp$ factors through $\D{D}_*^+\C{C}$, that is,
$$jp(c_1)=0\quad\textrm{ if }\quad c_1=[B\into A+B\onto A].$$
Let $c_0=[A+B]$ so that $c_0+\partial(c_1)=[A]+[B]$.
Then by
Lemma \ref{-c} (b) we have
\begin{align*}
  jp(c_1)&=-\alpha(c_0)+c_1+\alpha(c_0+\partial (c_1))\\&
={}-\alpha([A+B])\;+\;[B\into A+B\onto A]\;+\;\alpha([A]+[B])\,.
\end{align*}
This is zero by \eqref{alphaAB}.
\end{proof}

We finish this section with four lemmas which show useful relations in $\D{D}_*^+\C{C}$.


\begin{lem}\label{wcs+}
The following equality holds in $\D{D}_*^+\C{C}$.
\begin{eqnarray*}
&&[A+A'\;\into\; B+B'\;\onto\;
B/A+B'/A'\;\st{\sim}\l\;
C+C']\\&&\qquad=\;[A\into B\onto
B/A\st{\sim}\l C]^{[B']} +
[A'\into B'\onto B'/A'\st{\sim}\l C']
+\grupo{[A],[C']}.
\end{eqnarray*}
\end{lem}
\begin{proof}
Use \reln4, \reln{10} and Proposition \ref{weaknena} applied to the commutative diagram
$$\xymatrix{A'\ar@{>->}[r]\ar@{>->}[d]&B'\ar@{->>}[r]\ar@{>->}[d]&B'/A'\ar@{>->}[d]
&C'\ar[l]_\sim\ar@{>->}[d]
\\
A+A'\ar@{>->}[r]\ar@{->>}[d]&B+B'\ar@{->>}[r]\ar@{->>}[d]&B/A+B'/A'\ar@{->>}[d]
&C+C'\ar[l]_-\sim\ar@{->>}[d]\\
A\ar@{>->}[r]&B\ar@{->>}[r]&B/A&C\ar[l]_\sim
\\
A\ar@{=}[u]\ar@{>->}[r]&B\ar@{=}[u]\ar@{->>}[r]&B/A\ar@{=}[u]&C\ar@{=}[u]\ar[l]_\sim
}
$$
\end{proof}

As special cases we have
\begin{lem}\label{cs+}
The following equality holds in $\D{D}_*^+\C{C}$.
\begin{eqnarray*}
[A+A'\into B+B'\onto B/A+B'/A'] &&\\&&
\makebox[-5cm]{}
=\;\;[A\into B\onto B/A]^{[B']}
+[A'\into B'\onto B'/A']
+\grupo{[A],[B'/A']}.
\end{eqnarray*}
\end{lem}
\begin{lem}\label{we+}
The following equality holds in $\D{D}_*^+\C{C}$.
\begin{eqnarray*}
[A+ B\st{\sim}\r A'+ B']&=&
[A\st{\sim}\r A']^{[B']}+[B\st{\sim}\r B'].
\end{eqnarray*}
\end{lem}

We generalize \reln{$9'$} in the following lemma.

\begin{lem}\label{permu}
Let $\C{C}$ be a Waldhausen category with a functorial coproduct and let $A_1,\dots, A_n$ be objects in $\C{C}$.
Given a permutation $\sigma\in\sym{n}$ of $n$ elements we denote by
$$\sigma_{A_1,\ldots,A_n}\colon A_{\sigma_1}+\cdots + A_{\sigma_n}\st{\cong}\To A_1+\cdots+A_n$$
the isomorphism permuting the factors of the coproduct. Then the
following formula holds in $\D{D}_*^+\C{C}$.
\begin{eqnarray*}
[\sigma_{A_1,\ldots,A_n}]&=&\sum_{\substack{
i>j\\\sigma_i<\sigma_j}}\grupo{[A_{\sigma_i}],[A_{\sigma_j}]}.
\end{eqnarray*}
\end{lem}

\begin{proof}
The result holds for $n=1$ by $\reln4$. Suppose
$\sigma\in\sym{n}$ for $n\ge2$, and note that the isomorphism
$\sigma_{A_1,\ldots,A_n}$ factors naturally as
\begin{align*}
(\sigma'_{A_1,\ldots,A_{n-1}}+1)(1+\tau_{A_n,B})\colon
A_{\sigma_1}+\cdots + A_{\sigma_n} &\To A_{\sigma'_1}+\cdots +
A_{\sigma'_{n-1}}+A_n
\\& \To A_1+\cdots+A_{n-1}\,+\,A_n
\end{align*}
where
$(\sigma'_1,\ldots,\sigma'_{n-1},n)=(\sigma_1,\ldots,\widehat{\sigma_k},\ldots,\sigma_n,\sigma_k)$
and $B=A_{\sigma_{k+1}}+\cdots +A_{\sigma_n}$. Therefore by \reln4,
\reln6 and Lemma \ref{we+},
\begin{align*}
{}[\sigma_{A_1,\ldots,A_n}]&=[\sigma'_{A_1,\ldots,A_{n-1}}+1]+[1+\tau_{A_n,B}]
\\&=([\sigma'_{A_1,\ldots,A_{n-1}}]^{[A_n]}+0)+(0+[\tau_{A_n,B}]).
\end{align*}
By induction, \reln{$9'$} and Remark \ref{A+B} this is equal to
$$
\sum_{\substack{
p>q\\\sigma'_p<\sigma'_q}}\grupo{[A_{\sigma'_p}],[A_{\sigma'_q}]}
+ \grupo{[B],[A_n]} = \sum_{\substack{ i>j\\\sigma_i<\sigma_j<
n}}\grupo{[A_{\sigma_i}],[A_{\sigma_j}]} +\sum_{\substack{
i>j\\\sigma_i<\sigma_j=n}}\grupo{[A_{\sigma_i}],[A_{\sigma_j}]}$$
as required.
\end{proof}

\section{Proof of Theorem \ref{main}}

The morphism of stable quadratic modules $\D{D}_*F\colon\D{D}_*\C{C}\r\D{D}_*\C{D}$ induced by an exact functor $F\colon \C{C}\r\C{D}$, see \cite{1tK}, takes pairs of weak cofiber sequences to pairs of weak cofiber sequences,
\begin{eqnarray*}
(\D{D}_*F)\set{\begin{array}{c}\xymatrix@R=5pt@C=10pt{&&C_{1}\ar@{<-}[rd]^\sim&\\A\ar@<.7ex>@{>->}[r]\ar@<-.7ex>@{>->}[r]&B\ar@{->>}[ru]\ar@{->>}[rd]&&C\\
&&C_{2}\ar@{<-}[ru]_\sim&}\end{array}}
&=&
\set{\begin{array}{c}\xymatrix@R=5pt@C=10pt{&&F(C_{1})\ar@{<-}[rd]^\sim&\\F(A)\ar@<.7ex>@{>->}[r]\ar@<-.7ex>@{>->}[r]&F(B)\ar@{->>}[ru]\ar@{->>}[rd]&&F(C)\\
&&F(C_{2})\ar@{<-}[ru]_\sim&}\end{array}}.
\end{eqnarray*}
By \cite[Theorem 3.2]{1tK} exact equivalences of Waldhausen categories induce homotopy equivalences of stable quadratic modules, and hence isomorphisms in $K_1$. Therefore if the theorem holds for a certain Waldhausen category then it holds for all equivalent ones. In particular by Proposition \ref{hay} it is enough to prove the theorem for any Waldhausen category $\C{C}$ with a functorial coproduct such that the monoid of objects is freely generated by a set $\mathbb{S}$. By Theorem \ref{sumnor} we can work in the sum-normalized construction $\D{D}_*^+\C{C}$ in this case.

Any element $x\in \D{D}_1^+\C{C}$ is a sum of weak equivalences and cofiber sequences with coefficients $+1$ or $-1$. By Lemma \ref{we+} modulo the image of $\grupo{\cdot,\cdot}$ we can collect on the one hand all weak equivalences with coefficient $+1$ and on the other all weak equivalences with coefficient $-1$. Moreover, by Lemma  \ref{cs+} we can do the same for cofiber sequences. Therefore  the following equation holds modulo the image of $\grupo{\cdot,\cdot}$.
\begin{eqnarray*}
x&=&-[V_1\st{\sim}\r V_2]-[X_1\into X_2\onto X_3]\\
&&+[Y_1\into Y_2\onto Y_3]+[W_1\st{\sim}\r W_2]\\
\mbox{\scriptsize\reln4, \reln{10}, \ref{we+}, \ref{cs+}}\quad&=&
-[V_1+W_1\st{\sim}\r V_2+W_1]\\&&-[X_1+Y_1\into X_2+Y_3+Y_1\onto X_3+Y_3]\\
&&+[X_1+Y_1\into X_3+X_1+Y_2\onto X_3+Y_3]\\&&+[V_1+W_1\st{\sim}\r V_1+W_2]\\
\mbox{\scriptsize renaming}\quad&=&-[L\st{\sim}\r L_{1}]-[A\into E_{1}\onto D]\\
&&+[A\into E_{2}\onto D]+[L\st{\sim}\r L_{2}]\mod \grupo{\cdot,\cdot}.
\end{eqnarray*}

Suppose that $\partial(x)=0$ modulo commutators. Then modulo $[\cdot,\cdot]$
\begin{eqnarray}
\nonumber 0&=&-[L]+[L_{1}]-[A]-[D]+[E_{1}]\\&&-[E_{2}]+[D]+[A]-[L_{2}]+[L]\nonumber\\
\nonumber &=&[L_{1}]+[E_{1}]-[E_{2}]-[L_{2}]\mod[\cdot,\cdot],\\
\label{eso} \mbox{i.e. }[L_{1}+E_{1}]&=&[L_{2}+E_{2}]\mod[\cdot,\cdot].
\end{eqnarray}
The quotient of $\D{D}_0^+\C{C}$ by the commutator subgroup is the free abelian group on $\mathbb{S}$, hence by \eqref{eso}
there are objects $S_1,\dots,S_n\in\mathbb{S}$ and a permutation $\sigma\in\sym{n}$ such that
\begin{eqnarray*}
L_{1}+E_{1}&=&S_1+\cdots + S_n,\\
L_{2}+E_{2}&=&S_{\sigma_1}+\cdots + S_{\sigma_n},
\end{eqnarray*}
so there is an isomorphism
\begin{equation*}
\sigma_{S_1,\dots,S_n}\colon L_{2}+E_{2}\st{\sim}\To L_{1}+E_{1}.
\end{equation*}
Again by Lemmas \ref{we+} and \ref{cs+} modulo the image of $\grupo{\cdot,\cdot}$
\begin{eqnarray*}
x&=&-[L+D\st{\sim}\r L_{1}+D]-[A\into L_{1}+E_{1}\onto L_{1}+D]\\
&&+[A\into L_{2}+E_{2}\onto L_{2}+D]+[L+D\st{\sim}\r L_{2}+D]\\
\mbox{\scriptsize\ref{permu}}\quad&=&-[L+D\st{\sim}\r L_{1}+D]-[A\into L_{1}+E_{1}\onto L_{1}+D]\\
&&+[\sigma_{S_1,\dots,S_n}\colon L_{2}+E_{2}\st{\sim}\To L_{1}+E_{1}]\\
&&+[A\into L_{2}+E_{2}\onto L_{2}+D]+[L+D\st{\sim}\r L_{2}+D]\\
\mbox{\scriptsize\reln7}\quad&=&-[L+D\st{\sim}\r L_{1}+D]-[A\into L_{1}+E_{1}\onto L_{1}+D]\\
&&+[A\into L_{1}+E_{1}\onto L_{2}+D]+[L+D\st{\sim}\r L_{2}+D]\\
\mbox{\scriptsize renaming}\quad&=&-[C\st{\sim}\r C_{1}]-[A\into B\onto C_{1}]\\
&&+[A\into B\onto C_{2}]+[C\st{\sim}\r C_{2}]\\
&=&-[C\st{\sim}\r C_{1}]^{[A]}-[A\into B\onto C_{1}]\\
&&+[A\into B\onto C_{2}]+[C\st{\sim}\r C_{2}]^{[A]}\\
&=&-[A\into B\onto C_{1}\st{\sim}\l C]+[A\into B\onto C_{2}\st{\sim}\l C]\mod \grupo{\cdot,\cdot},
\end{eqnarray*}
i.e. there is $y\in\D{D}_1^+\C{C}$ in the image of $\grupo{\cdot,\cdot}$ such that
\begin{eqnarray*}
x&=&\set{\begin{array}{c}\xymatrix@R=5pt@C=10pt{&&C_{1}\ar@{<-}[rd]^\sim&\\
A\ar@<.7ex>@{>->}[r]\ar@<-.7ex>@{>->}[r]&B\ar@{->>}[ru]\ar@{->>}[rd]&&C\\
&&C_{2}\ar@{<-}[ru]_\sim&}\end{array}}+y.
\end{eqnarray*}

Now assume that $\partial(x)=0$. Since the pair of weak cofiber sequences is also in the kernel of $\partial$ we have $\partial(y)=0$.
In order to give the next step we need a technical lemma.

\begin{lem}\label{kita}
Let $C_*$ be a stable quadratic module such that $C_0$ is a free group of nilpotency class $2$. Then
any element $y\in \ker\partial\cap\operatorname{Image}\grupo{\cdot,\cdot}$ is of the form $y=\grupo{a,a}$ for some $a\in C_0$.
\end{lem}

\begin{proof}
For any abelian group $A$ let $\hat{\otimes}^2A$ be the quotient of the tensor square $A\otimes A$ by the relations $a\otimes b+b\otimes a=0$, $a,b\in A$, and let $\wedge^2A$ be the quotient of $A\otimes A$ by the relations $a\otimes a=0$, $a\in A$, which is also a quotient of $\hat{\otimes}^2A$. The projection of $a\otimes b\in A\otimes A$ to $\hat{\otimes}^2A$ and $\wedge^2A$ is denoted by $a\hat{\otimes}b$ and $a\wedge b$, respectively.

There is a commutative diagram of group homomorphisms
$$\xymatrix{&C_0^\abb\otimes C_0^\abb\ar@{->>}[d]\ar@/^15pt/@{->>}[rd]\ar@/_25pt/@{-->}[dd]|{\quad}_<(.7){\grupo{\cdot,\cdot}}&\\
C_0^\abb\otimes\Z/2\;\ar@{^{(}->}[r]^(0.52){\bar{\tau}}&\hat{\otimes}^2C_0^\abb\ar@{->>}[r]^-q
\ar[d]^{c_1}&\wedge^2C_0^\abb\ar[d]^{c_0}\\
&C_1\ar[r]^-\partial &C_0}$$
where 
$\bar{\tau}(a\otimes1)=a\hat{\otimes} a$, the factorization $c_1$ of $\grupo{\cdot,\cdot}$ is given by $c_1(a\hat\otimes b)=\grupo{a,b}$, which is well defined by Definition \ref{ob} (3), and $c_0(a\wedge b)=[b,a]$ is well known to be injective in the case $C_0$ is free of nilpotency class $2$.

Moreover, $C^\abb_0$ is a free abelian group, hence the middle row is a short exact sequence, see \cite[I.4]{ch4c}. Therefore any element $y\in C_1$ which is both in the image of $c_1$ and the kernel of $\partial$ is in the image of $c_1\bar\tau$ as required.
\end{proof}

The group $\D{D}_0^+\C{C}$ is free of nilpotency class $2$
, therefore by the previous lemma $y=\grupo{a,a}$ for some $a\in \D{D}_0^+\C{C}$. By the laws of a stable quadratic module the element $\grupo{a,a}$ only depends on $a$ mod $2$, compare \cite[Definition 1.6]{1tK}, and so we can suppose that $a$ is a sum of basis elements with coefficient $+1$. Hence by Remark \ref{A+B} we can take $a=[M]$ for some object $M$ in $\C{C}$,
\begin{eqnarray*}
y&=&\grupo{[M],[M]}.
\end{eqnarray*}
The element $\grupo{[M],[M]}$ is itself a pair of weak cofiber sequences
\begin{eqnarray*}
\grupo{[M],[M]}&=&\left\{\begin{array}{c}\xymatrix@R=5pt@C=10pt{&&M\ar@{<-}[rd]^\sim_1&\\
M\ar@<.7ex>@{>->}[r]^-{i_2}\ar@<-.7ex>@{>->}[r]_-{i_1}&M+M\ar@{->>}[ru]^-{p_1}\ar@{->>}[rd]_-{p_2}&&M\\
&&M\ar@{<-}[ru]_\sim^1&}\end{array}\right\},
\end{eqnarray*}
therefore $x$ is also a pair of weak cofiber sequences,  by Proposition \ref{spwcs}, and the proof of Theorem \ref{main} is complete.

\section{Comparison with Nenashev's approach for exact categories}

For any Waldhausen category $\C{C}$ we denote by $K^\mathrm{wcs}_1\,\C{C}$ the abelian group generated by pairs of weak cofiber sequences
\begin{equation*}
\left\{
\begin{array}{c}
\xymatrix@R=5pt@C=10pt{&&C_{1}\ar@{<-}[rd]^\sim&\\
A\ar@<.7ex>@{>->}[r]\ar@<-.7ex>@{>->}[r]&B\ar@{->>}[ru]\ar@{->>}[rd]&&C\\
&&C_{2}\ar@{<-}[ru]_\sim&}
\end{array}
\right\}
\end{equation*}
modulo the relations \rels1 and \rels2 in Section \ref{relations}. Theorem \ref{larel} and Proposition \ref{las0} show the existence of a natural homomorphism
\begin{equation}\label{sur}
K^\mathrm{wcs}_1\,\C{C}\onto K_1\C{C}
\end{equation}
which is surjective by Theorem \ref{main}.

Given an exact category $\C{E}$ Nenashev defines in \cite{dsesk1} an abelian group $D(\C{E})$ by generators and relations which surjects naturally to $K_1\C{E}$. Moreover, he shows in \cite{K1gr} that this natural surjection is indeed a natural isomorphism
\begin{eqnarray}\label{neniso}
D(\C{E})&\cong&K_1\C{E}.
\end{eqnarray}
Generators of $D(\C{E})$ are pairs of short exact sequences
\begin{eqnarray*}
\{\xymatrix@R=5pt@C=10pt{
A\ar@<.7ex>@{>->}[r]^{j_1}\ar@<-.7ex>@{>->}[r]_{j_2}&B\ar@<.7ex>@{->>}[r]^{r_1}\ar@<-.7ex>@{->>}[r]_{r_2}&C}\}.
\end{eqnarray*}
They satisfy two kind of relations. The first, analogous to \rels2, says that a generator vanishes provided $j_1=j_2$ and $r_1=r_2$. The second is a simplification of \rels1: given six pairs of short exact sequences
\begin{equation*}
\begin{array}{ccc}
\xymatrix@R=5pt@C=10pt{
A'\ar@<.7ex>@{>->}[r]^{j^A_1}\ar@<-.7ex>@{>->}[r]_{j^A_2}&A\ar@{->>}@<.7ex>[r]^{r^A_1}\ar@{->>}@<-.7ex>[r]_{r^A_2}&A'',}
&
\xymatrix@R=5pt@C=10pt{
B'\ar@<.7ex>@{>->}[r]^{j^B_1}\ar@<-.7ex>@{>->}[r]_{j^B_2}&B\ar@{->>}@<.7ex>[r]^{r^B_1}\ar@{->>}@<-.7ex>[r]_{r^B_2}&B'',}
&
\xymatrix@R=5pt@C=10pt{
C'\ar@<.7ex>@{>->}[r]^{j^C_1}\ar@<-.7ex>@{>->}[r]_{j^C_2}&C\ar@{->>}@<.7ex>[r]^{r^C_1}\ar@{->>}@<-.7ex>[r]_{r^C_2}&C'',}
\\
\xymatrix@R=5pt@C=10pt{
A'\ar@<.7ex>@{>->}[r]^{j'_1}\ar@<-.7ex>@{>->}[r]_{j'_2}&B'\ar@{->>}@<.7ex>[r]^{r'_1}\ar@{->>}@<-.7ex>[r]_{r'_2}&C',}
&
\xymatrix@R=5pt@C=10pt{
A\ar@<.7ex>@{>->}[r]^{j_1}\ar@<-.7ex>@{>->}[r]_{j_2}&B\ar@{->>}@<.7ex>[r]^{r_1}\ar@{->>}@<-.7ex>[r]_{r_2}&C,}
&
\xymatrix@R=5pt@C=10pt{
A''\ar@<.7ex>@{>->}[r]^{j''_1}\ar@<-.7ex>@{>->}[r]_{j''_2}&B''\ar@{->>}@<.7ex>[r]^{r''_1}\ar@{->>}@<-.7ex>[r]_{r''_2}&C'',}
\end{array}
\end{equation*}
denoted for simplicity as
$\lambda^A$, $\lambda^B$, $\lambda^C$, $\lambda'$, $\lambda$, and $\lambda''$, such that the diagram
\begin{equation*}
\xymatrix{
A'\ar@{>->}[r]^-{j^A_i}\ar@{>->}[d]_-{j'_i}&
A\ar@{>->}[d]_-{j_i}\ar@{->>}[r]^-{r^A_{i}}&
A''\ar@{>->}[d]^-{j''_{i}}
\\
B'\ar@{>->}[r]^-{j^B_i}\ar@{->>}[d]_-{r'_{i}}&
B\ar@{->>}[r]^-{r^B_{i}}\ar@{->>}[d]_-{r_{i}}&
B'i\ar@{->>}[d]^-{r''_{i}}
\\
C'\ar@{>->}[r]_-{j^C_i}&C
\ar@{->>}[r]_-{r^C_i}&C''}
\end{equation*}
commutes for $i=1,2$ then
\begin{eqnarray*}
\set{\lambda^A}-\set{\lambda^B}+\set{\lambda^C}&=&\set{\lambda'}-\set{\lambda}+\set{\lambda''}.
\end{eqnarray*}

\begin{prop}\label{comp}
For any exact category $\C{E}$ there is a natural isomorphism
\begin{eqnarray*}
D(\C{E})&\cong&K_1^\mathrm{wcs}\;\C{E}.
\end{eqnarray*}
\end{prop}

\begin{proof}
There is a clear homomorphism $D(\C{E})\r K_1^\mathrm{wcs}\;\C{E}$ defined by
\begin{eqnarray*}
\{\xymatrix@R=5pt@C=10pt{
A\ar@<.7ex>@{>->}[r]^{j_1}\ar@<-.7ex>@{>->}[r]_{j_2}&B\ar@<.7ex>@{->>}[r]^{r_1}\ar@<-.7ex>@{->>}[r]_{r_2}&C}\}
&\mapsto&
\left\{
\begin{array}{c}
\xymatrix@R=5pt@C=10pt{&&C\ar@{<-}[rd]^\sim_1&\\
A\ar@<.7ex>@{>->}[r]^{j_1}\ar@<-.7ex>@{>->}[r]_{j_2}&B\ar@{->>}[ru]^{r_1}\ar@{->>}[rd]_{r_2}&&C\\
&&C\ar@{<-}[ru]^1_\sim&}
\end{array}
\right\}
\end{eqnarray*}
Exact categories regarded as Waldhausen categories have the particular feature that all weak equivalences are isomorphisms, hence one can also define a homomorphism $K_1^\mathrm{wcs}\,\C{E}\r D(\C{E})$ as
\begin{eqnarray*}
\left\{
\begin{array}{c}
\xymatrix@R=5pt@C=10pt{&&C_1\ar@{<-}[rd]^\cong_{w_1}&\\
A\ar@<.7ex>@{>->}[r]^{j_1}\ar@<-.7ex>@{>->}[r]_{j_2}&B\ar@{->>}[ru]^{r_1}\ar@{->>}[rd]_{r_2}&&C\\
&&C_2\ar@{<-}[ru]^{w_2}_\cong&}\end{array}
\right\}
&\mapsto&
\{\xymatrix@R=5pt@C=10pt{
A\ar@<.7ex>@{>->}[r]^{j_1}\ar@<-.7ex>@{>->}[r]_{j_2}&B\ar@<.7ex>@{->>}[rr]^{w_1^{-1}r_1}\ar@<-.7ex>@{->>}[rr]_{w_2^{-1}r_2}&&C}\}.
\end{eqnarray*}
We leave the reader to check the compatibiliy with the relations and the fact that these two homomorphisms are inverse of each other.
\end{proof}

\begin{rem}
The composite of the isomorphism in Proposition \ref{comp} with the natural epimorphism (\ref{sur}) for $\C{C}=\C{E}$ coincides with Nenashev's isomorphism (\ref{neniso}), therefore the natural epimorphism (\ref{sur}) is an isomorphism when $\C{C}=\C{E}$ is an exact category.
\end{rem}

\section{Weak cofiber sequences and the stable Hopf map}

After the previous section one could conjecture that the natural epimorphism
\begin{equation*}
K^\mathrm{wcs}_1\; \C{C}\onto K_1\C{C}
\end{equation*}
in (\ref{sur}) is an isomorphism not only for exact categories but for any Waldhausen category $\C{C}$. In order to support this conjecture we show the following result.

\begin{thm}\label{factor}
For any Waldhausen category $\C{C}$ there is a natural homomorphism
$$\phi\colon K_0\C{C}\otimes\Z/2\To K^\mathrm{wcs}_1\; \C{C}$$
which composed with (\ref{sur}) yields the homomorphism $\cdot\,\eta\colon K_0\C{C}\otimes\Z/2\r K_1 \C{C}$ determined by the
action of the stable Hopf map $\eta\in\pi_*^s$ in the stable homotopy groups of spheres.
\end{thm}

For the proof we use the following lemma.

\begin{lem}\label{compot}
The following relation holds in $K^\mathrm{wcs}_1\; \C{C}$
\begin{eqnarray*}
\left\{\begin{array}{c}\xymatrix@R=5pt@C=10pt{&&A''\ar@{<-}[rd]_\sim^{w'_1w_1}&\\
{*}\ar@<.7ex>@{>->}[r]\ar@<-.7ex>@{>->}[r]&A''\ar@{->>}[ru]^1\ar@{->>}[rd]_1&&A\\
&&A''\ar@{<-}[ru]^\sim_{w'_2w_2}&}\end{array}\right\}&=&
\left\{\begin{array}{c}\xymatrix@R=5pt@C=10pt{&&A'\ar@{<-}[rd]_\sim^{w_1}&\\
{*}\ar@<.7ex>@{>->}[r]\ar@<-.7ex>@{>->}[r]&A'\ar@{->>}[ru]^1\ar@{->>}[rd]_1&&A\\
&&A'\ar@{<-}[ru]^\sim_{w_2}&}\end{array}\right\}                    \\&&
+\left\{\begin{array}{c}\xymatrix@R=5pt@C=10pt{&&A''\ar@{<-}[rd]_\sim^{w'_1}&\\
{*}\ar@<.7ex>@{>->}[r]\ar@<-.7ex>@{>->}[r]&A''\ar@{->>}[ru]^1\ar@{->>}[rd]_1&&A'\\
&&A''\ar@{<-}[ru]^\sim_{w'_2}&}\end{array}\right\}.
\end{eqnarray*}
\end{lem}

\begin{proof}
Use relations \rels1 and \rels2 applied to the following pairs of weak cofiber sequences
\begin{equation*}
\begin{array}{ccc}
\begin{array}{c}\xymatrix@R=5pt@C=10pt{&&{*}\ar@{<-}[rd]_\sim&\\
{*}\ar@<.7ex>@{>->}[r]\ar@<-.7ex>@{>->}[r]&{*}\ar@{->>}[ru]\ar@{->>}[rd]&&{*},\\
&&{*}\ar@{<-}[ru]^\sim&}\end{array}
&
\begin{array}{c}\xymatrix@R=5pt@C=10pt{&&A''\ar@{<-}[rd]_\sim^{w'_1w_1}&\\
{*}\ar@<.7ex>@{>->}[r]\ar@<-.7ex>@{>->}[r]&A''\ar@{->>}[ru]^1\ar@{->>}[rd]_1&&A,\\
&&A''\ar@{<-}[ru]^\sim_{w'_2w_2}&}\end{array}
&
\begin{array}{c}\xymatrix@R=5pt@C=10pt{&&A'\ar@{<-}[rd]_\sim^{w_1}&\\
{*}\ar@<.7ex>@{>->}[r]\ar@<-.7ex>@{>->}[r]&A'\ar@{->>}[ru]^1\ar@{->>}[rd]_1&&A,\\
&&A'\ar@{<-}[ru]^\sim_{w_2}&}\end{array}
\\
\begin{array}{c}\xymatrix@R=5pt@C=10pt{&&{*}\ar@{<-}[rd]_\sim&\\
{*}\ar@<.7ex>@{>->}[r]\ar@<-.7ex>@{>->}[r]&{*}\ar@{->>}[ru]\ar@{->>}[rd]&&{*},\\
&&{*}\ar@{<-}[ru]^\sim&}\end{array}
&
\begin{array}{c}\xymatrix@R=5pt@C=10pt{&&A''\ar@{<-}[rd]_\sim^{w'_1}&\\
{*}\ar@<.7ex>@{>->}[r]\ar@<-.7ex>@{>->}[r]&A''\ar@{->>}[ru]^1\ar@{->>}[rd]_1&&A',\\
&&A''\ar@{<-}[ru]^\sim_{w'_2}&}\end{array}
&
\begin{array}{c}\xymatrix@R=5pt@C=10pt{&&A\ar@{<-}[rd]_\sim^{1}&\\
{*}\ar@<.7ex>@{>->}[r]\ar@<-.7ex>@{>->}[r]&A\ar@{->>}[ru]^1\ar@{->>}[rd]_1&&A.\\
&&A\ar@{<-}[ru]^\sim_{1}&}\end{array}
\end{array}
\end{equation*}
\end{proof}

Now we are ready to prove Theorem \ref{factor}.

\begin{proof}[Proof of Theorem \ref{factor}]
We define the homomorphism by
\begin{equation*}
\phi[A]=\left\{
\!\!\!\!
\begin{array}{c}\xymatrix@R=5pt@C=10pt{&&A\ar@{<-}[rd]_\sim^1&\\
A\ar@<.7ex>@{>->}[r]^-{i_2}\ar@<-.7ex>@{>->}[r]_-{i_1}&A\vee A\ar@{->>}[ru]^{p_1}\ar@{->>}[rd]_{p_2}&&A\\
&&A\ar@{<-}[ru]^\sim_1&}\end{array}
\!\!\!\!
\right\}   
                       =
\left\{
\!\!\!\!
\begin{array}{c}\xymatrix@R=5pt@C=10pt{&&A\vee A\ar@{<-}[rd]_\sim^1&\\
{*}\ar@<.7ex>@{>->}[r]\ar@<-.7ex>@{>->}[r]&A\vee A\ar@{->>}[ru]^{1}\ar@{->>}[rd]_{1}&&A\vee A\\
&&A\vee A\ar@{<-}[ru]^\sim_{\;\tau_{A,A}}&}\end{array}
\!\!\!\!
\right\}.
\end{equation*}
Here the second equality follows from
\rels1 and \rels2
applied to the following pairs of weak cofiber sequences
\begin{equation*}
\begin{array}{ccc}
\!\!\!\!
\begin{array}{c}\xymatrix@R=5pt@C=10pt{&&A\ar@{<-}[rd]_\sim^{1}&\\
{*}\ar@<.7ex>@{>->}[r]\ar@<-.7ex>@{>->}[r]&A\ar@{->>}[ru]^1\ar@{->>}[rd]_1&&A,\\
&&A\ar@{<-}[ru]^\sim_{1}&}\end{array}
&
\!\!\!\!
\begin{array}{c}\xymatrix@R=5pt@C=10pt{&&A\vee A\ar@{<-}[rd]_\sim^1&\\
{*}\ar@<.7ex>@{>->}[r]\ar@<-.7ex>@{>->}[r]&A\vee A\ar@{->>}[ru]^{1}\ar@{->>}[rd]_{1}&&A\vee A,\\
&&A\vee A\ar@{<-}[ru]^\sim_{\;\tau_{A,A}}&}\end{array}
&
\!\!\!\!\!\!\!\!
\begin{array}{c}\xymatrix@R=5pt@C=10pt{&&A\ar@{<-}[rd]_\sim^{1}&\\
{*}\ar@<.7ex>@{>->}[r]\ar@<-.7ex>@{>->}[r]&A\ar@{->>}[ru]^1\ar@{->>}[rd]_1&&A,\\
&&A\ar@{<-}[ru]^\sim_{1}&}\end{array}
\\
\!\!\!\!\!\!\!\!
\begin{array}{c}\xymatrix@R=5pt@C=10pt{&&{*}\ar@{<-}[rd]_\sim&\\
{*}\ar@<.7ex>@{>->}[r]\ar@<-.7ex>@{>->}[r]&{*}\ar@{->>}[ru]\ar@{->>}[rd]&&{*},\\
&&{*}\ar@{<-}[ru]^\sim&}\end{array}
&
\!\!\!\!
\begin{array}{c}\xymatrix@R=5pt@C=10pt{&&A\ar@{<-}[rd]_\sim^1&\\
A\ar@<.7ex>@{>->}[r]^-{i_2}\ar@<-.7ex>@{>->}[r]_-{i_1}&A\vee A\ar@{->>}[ru]^{p_1}\ar@{->>}[rd]_{p_2}&&A,\\
&&A\ar@{<-}[ru]^\sim_1&}\end{array}
&
\!\!\!\!
\begin{array}{c}\xymatrix@R=5pt@C=10pt{&&A\ar@{<-}[rd]_\sim^1&\\
A\ar@<.7ex>@{>->}[r]^-{i_2}\ar@<-.7ex>@{>->}[r]_-{i_2}&A\vee A\ar@{->>}[ru]^{p_1}\ar@{->>}[rd]_{p_1}&&A.\\
&&A\ar@{<-}[ru]^\sim_1&}\end{array}
\end{array}
\end{equation*}

Given a cofiber sequence $A\st{j}\into B\st{r}\onto C$ the equation $\phi[C]+\phi[A]=\phi[B]$ follows from
\rels1 and \rels2
applied to the following pairs of weak cofiber sequences
\begin{equation*}
\begin{array}{ccc}
\!\!\!\!
\begin{array}{c}\xymatrix@R=5pt@C=10pt{&&A\ar@{<-}[rd]_\sim^1&\\
A\ar@<.7ex>@{>->}[r]^-{i_2}\ar@<-.7ex>@{>->}[r]_-{i_1}&A\vee A\ar@{->>}[ru]^{p_1}\ar@{->>}[rd]_{p_2}&&A,\\
&&A\ar@{<-}[ru]^\sim_1&}\end{array}
&
\!\!\!\!\!\!\!\!\!\!\!\!\!\!\!\!\!\!\!
\begin{array}{c}\xymatrix@R=5pt@C=10pt{&&B\ar@{<-}[rd]_\sim^1&\\
B\ar@<.7ex>@{>->}[r]^-{i_2}\ar@<-.7ex>@{>->}[r]_-{i_1}&B\vee B\ar@{->>}[ru]^{p_1}\ar@{->>}[rd]_{p_2}&&B,\\
&&B\ar@{<-}[ru]^\sim_1&}\end{array}
&
\!\!\!\!\!\!\!\!\!\!\!\!\!\!\!\!\!\!\!\!\!\!\!
\begin{array}{c}\xymatrix@R=5pt@C=10pt{&&C\ar@{<-}[rd]_\sim^1&\\
C\ar@<.7ex>@{>->}[r]^-{i_2}\ar@<-.7ex>@{>->}[r]_-{i_1}&C\vee C\ar@{->>}[ru]^{p_1}\ar@{->>}[rd]_{p_2}&&C,\\
&&C\ar@{<-}[ru]^\sim_1&}\end{array}
\\
\!\!\!\!\!\!\!\!
\begin{array}{c}\xymatrix@R=5pt@C=10pt{&&C\ar@{<-}[rd]_\sim^1&\\
A\ar@<.7ex>@{>->}[r]^j\ar@<-.7ex>@{>->}[r]_j&B\ar@{->>}[ru]^r\ar@{->>}[rd]_r&&C,\\
&&C\ar@{<-}[ru]^\sim_1&}\end{array}
&
\!\!\!\!\!\!\!\!\!\!\!\!\!\!\!\!
\begin{array}{c}\xymatrix@R=5pt@C=10pt{&&C\vee C\ar@{<-}[rd]_\sim^1&\\
A\vee A\ar@<.7ex>@{>->}[r]^{j\vee j}\ar@<-.7ex>@{>->}[r]_{j\vee j}&B\vee B\ar@{->>}[ru]^{r\vee r}\ar@{->>}[rd]_{r\vee r}&&C\vee C,\\
&&C\vee C\ar@{<-}[ru]^\sim_1&}\end{array}
&
\!\!\!\!\!\!\!\!\!\!\!\!
\begin{array}{c}\xymatrix@R=5pt@C=10pt{&&C\ar@{<-}[rd]_\sim^1&\\
A\ar@<.7ex>@{>->}[r]^j\ar@<-.7ex>@{>->}[r]_j&B\ar@{->>}[ru]^r\ar@{->>}[rd]_r&&C.\\
&&C\ar@{<-}[ru]^\sim_1&}\end{array}
\end{array}
\end{equation*}

Given a weak equivalence $w\colon A\r A'$ the equation $\phi[A]=\phi[A']$ follows from
\rels1 and \rels2
applied to the following pairs of weak cofiber sequences
\begin{equation*}
\begin{array}{ccc}
\!\!\!\!
\begin{array}{c}\xymatrix@R=5pt@C=10pt{&&{*}\ar@{<-}[rd]_\sim&\\
{*}\ar@<.7ex>@{>->}[r]\ar@<-.7ex>@{>->}[r]&{*}\ar@{->>}[ru]\ar@{->>}[rd]&&{*},\\
&&{*}\ar@{<-}[ru]^\sim&}\end{array}
&
\!\!\!\!
\begin{array}{c}\xymatrix@R=5pt@C=10pt{&&A'\ar@{<-}[rd]_\sim^1&\\
A'\ar@<.7ex>@{>->}[r]^-{i_2}\ar@<-.7ex>@{>->}[r]_-{i_1}&A'\vee A'\ar@{->>}[ru]^{p_1}\ar@{->>}[rd]_{p_2}&&A',\\
&&A'\ar@{<-}[ru]^\sim_1&}\end{array}
&
\!\!\!\!\!\!\!\!\!\!\!
\begin{array}{c}\xymatrix@R=5pt@C=10pt{&&A\ar@{<-}[rd]_\sim^1&\\
A\ar@<.7ex>@{>->}[r]^-{i_2}\ar@<-.7ex>@{>->}[r]_-{i_1}&A\vee A\ar@{->>}[ru]^{p_1}\ar@{->>}[rd]_{p_2}&&A,\\
&&A\ar@{<-}[ru]^\sim_1&}\end{array}
\\
\!\!\!\!\!\!\!
\begin{array}{c}\xymatrix@R=5pt@C=10pt{&&A'\ar@{<-}[rd]_\sim^{w}&\\
{*}\ar@<.7ex>@{>->}[r]\ar@<-.7ex>@{>->}[r]&A'\ar@{->>}[ru]^1\ar@{->>}[rd]_1&&A,\\
&&A'\ar@{<-}[ru]^\sim_{w}&}\end{array}
&
\!\!\!\!
\begin{array}{c}\xymatrix@R=5pt@C=10pt{&&A'\vee A'\ar@{<-}[rd]_\sim^{w\vee w}&\\
{*}\ar@<.7ex>@{>->}[r]\ar@<-.7ex>@{>->}[r]&A'\vee A'\ar@{->>}[ru]^1\ar@{->>}[rd]_1&&A\vee A,\\
&&A'\vee A'\ar@{<-}[ru]^\sim_{w\vee w}&}\end{array}
&
\!\!\!\!
\begin{array}{c}\xymatrix@R=5pt@C=10pt{&&A'\ar@{<-}[rd]_\sim^{w}&\\
{*}\ar@<.7ex>@{>->}[r]\ar@<-.7ex>@{>->}[r]&A'\ar@{->>}[ru]^1\ar@{->>}[rd]_1&&A.\\
&&A'\ar@{<-}[ru]^\sim_{w}&}\end{array}
\end{array}
\end{equation*}

The equation $2\phi[A]=0$ follows from
\rels2
and Lemma \ref{compot} by using the second formula for $\phi[A]$ above and the fact that $\tau_{A,A}^2=1$.

We have already shown that $\phi$ is a well-defined homomorphism. The composite of $\phi$ and (\ref{sur}) coincides with the action of the stable Hopf map by the main result of \cite{1tK}.
\end{proof}

\appendix

\section{Free stable quadratic modules and presentations}

Free stable quadratic modules, and also stable quadratic modules
defined by generators and relations, can be characterized up to
isomorphism by obvious universal properties. In this appendix we
give more explicit constructions of these notions.

Let $\C{squad}$ be the category of stable quadratic modules and
let
$$U\colon\C{squad}\To\C{Set}\times\C{Set}$$
be the forgetful functor, $U(C_*)=(C_0,C_1)$. The functor $U$ has
a left adjoint $F$, and a stable quadratic module $F(E_0,E_1)$ is
called a {\em free  stable quadratic module} on the sets $E_0$ and
$E_1$. In order to give an explicit description of $F(E_0,E_1)$ we
fix some notation. Given a set $E$ we denote by $\grupo{E}$ the
\emph{free group} with basis $E$, and by $\grupo{E}^{\abb}$ the
\emph{free abelian group} with basis $E$. The \emph{free group of
nilpotency class $2$} with basis $E$, denoted by
$\grupo{E}^{\ni}$, is the quotient of $\grupo{E}$ by triple
commutators. Given a pair of sets $E_0$ and $E_1$, we write
$E_0\cup
\partial E_1$ for the set whose elements are the symbols $e_0$ and
$\partial e_1$ for each $e_0\in E_0$, $e_1\in E_1$.

To define $F(E_0,E_1)$, consider the groups
\begin{eqnarray*}
F(E_0,E_1)_0&=&\grupo{E_0\cup \partial E_1}^\ni,\\
F(E_0,E_1)_1&=&(\grupo{E_0}^{\abb}\otimes\Z/2)\times\ker\delta.
\end{eqnarray*}
Here $\delta\colon F(E_0,E_1)_0\r \grupo{E_0}^{ab}$ is the
homomorphism given by $\delta e_0=e_0$ and $\delta\partial e_1=0$.
In the notation of the proof of Lemma \ref{kita} there are
isomorphisms
\begin{eqnarray*}
\ker\delta &\cong& \wedge^2\grupo{E_0}^{\abb}\times
\grupo{E_0\times
E_1}^\abb\times\grupo{E_1}^\ni,\\
(\grupo{E_0}^{\abb}\otimes\Z/2)\times\ker\delta &\cong&
\hat\otimes^2\grupo{E_0}^{\abb}\times \grupo{E_0\times
E_1}^\abb\times\grupo{E_1}^\ni,
\end{eqnarray*}
and intuitively we think of $F(E_0,E_1)_1$ as a group generated by
symbols $\grupo{e_0,e_0'}$, $\grupo{e_0,\partial e_1}$ and $e_1$.
The symbol $\grupo{\partial e_1,\partial e_1'}$ is unnecessary
since it will be given by the commutator $[e_1',e_1]$.

We define structure homomorphisms on $F(E_0,E_1)$ as follows. The
boundary
$$\partial\colon F(E_0,E_1)_1\To F(E_0,E_1)_0$$ is the projection
onto $\ker\delta$ followed by the inclusion of the kernel. The
bracket
$$\grupo{\cdot,\cdot}\colon F(E_0,E_1)_0^{ab}\otimes
F(E_0,E_1)_0^{ab}\To F(E_0,E_1)_1
$$ is
given by the product of the following two homomorphisms,
$$c'\colon F(E_0,E_1)_0^{ab}\otimes F(E_0,E_1)_0^{ab}\To \grupo{E_0}^{ab}\otimes\Z/2,$$
defined on the generators $x,y\in E_0\cup\partial E_1$ by
$c'(x,y)=x\otimes 1$ if $x=y\in E_0$ and $c'(x,y)=0$ otherwise,
and
$$c\colon F(E_0,E_1)_0^{ab}\otimes F(E_0,E_1)_0^{ab}\To \ker\delta$$
induced by the commutator bracket, $c(a,b)=[b,a]$.

It is now straightforward to define explicitly the stable
quadratic module $C_\ast$ presented by {\em generators} $E_i$ and
{\em relations} $R_i\subset F(E_0,E_1)_i$ in degrees $i=0,1$, by
\begin{eqnarray*}
C_0&=&F(E_0,E_1)_0/N_0,\\
C_1&=&F(E_0,E_1)_1/N_1.
\end{eqnarray*}
Here $N_0\subset F(E_0,E_1)_0$ is the normal subgroup generated by
the elements of $R_0$ and $\partial R_1$, and $N_1\subset
F(E_0,E_1)_1$ is the normal subgroup generated by the elements of
$R_1$ and $\grupo{F(E_0,E_1)_0,N_0}$. The boundary and bracket on
$F(E_0,E_1)$ induce a stable quadratic module structure on $C_*$
which satisfies the following universal property: given a stable
quadratic module $C_*'$, any 
pair of
functions $E_i\r C_i'$ $(i=0,1)$ such that the
induced morphism $F(E_0,E_1)\r C_*'$ annihilates $R_0$ and $R_1$
induces a unique morphism $C_*\r C_*'$.

In \cite{ch4c} Baues considers the \emph{totally free} stable
quadratic module $C_*$ with basis given by a function $g\colon
E_1\r\grupo{E_0}^{nil}$. In the language of this paper $C_*$ is
the stable quadratic module with generators $E_i$ in degree
$i=0,1$ and degree $0$ relations $\partial(e_1)=g(e_1)$ for all
$e_1\in E_1$.

\bibliographystyle{amsalpha}
\bibliography{Fernando}
\end{document}